\documentclass[english,12pt]{smfart}
\usepackage{amssymb}
\usepackage{amsfonts}
\usepackage{amscd}
\usepackage[all]{xypic}

\CompileMatrices

\swapnumbers

\def\ord{{\rm ord}}
\def\ac{{\overline{\rm ac}}}

\def\11{{\mathbf 1}}
\def\AA{{\mathbf A}}

\def\CC{{\mathbf C}}

\def\FF{{\mathbf F}}
\def\GG{{\mathbf G}}

\def\NN{{\mathbf N}}

\def\QQ{{\mathbf Q}}
\def\RR{{\mathbf R}}

\def\ZZ{{\mathbf Z}}

\def\cI{{\mathcal I}}

\def\cO{{\mathcal O}}

\mathchardef\alphag="7C0B \mathchardef\betag="7C0C
\mathchardef\gammag="7C0D \mathchardef\deltag="7C0E
\mathchardef\varepsilong="7C22 \mathchardef\varphig="7C27
\mathchardef\psig="7C20 \mathchardef\zetag="7C10
\mathchardef\epsilong="7C0F \mathchardef\rhog="7C1A
\mathchardef\taug="7C1C \mathchardef\upsilong="7C1D
\mathchardef\iotag="7C13 \mathchardef\thetag="7C12
\mathchardef\pig="7C19 \mathchardef\sigmag="7C1B
\mathchardef\etag="7C11 \mathchardef\omegag="7C21
\mathchardef\kappag="7C14 \mathchardef\lambdag="7C15
\mathchardef\mug="7C16 \mathchardef\xig="7C18
\mathchardef\chig="7C1F \mathchardef\nug="7C17
\mathchardef\varthetag="7C23 \mathchardef\varpig="7C24
\mathchardef\varrhog="7C25 \mathchardef\varsigmag="7C26
\mathchardef\Omegag="7C0A \mathchardef\Thetag="7C02
\mathchardef\Sigmag="7C06 \mathchardef\Deltag="7C01
\mathchardef\Phig="7C08 \mathchardef\Gammag="7C00
\mathchardef\Psig="7C09 \mathchardef\Lambdag="7C03
\mathchardef\Xig="7C04 \mathchardef\Pig="7C05
\mathchardef\Upsilong="7C07

\newtheorem{theorem}[subsection]{Theorem}
\newtheorem{lem}[subsection]{Lemma}
\newtheorem{cor}[subsection]{Corollary}
\newtheorem{prop}[subsection]{Proposition}

\theoremstyle{definition}

\newtheorem{def-prop}[subsection]{Proposition-Definition}
\newtheorem{def-theorem}[subsection]{Theorem-Definition}
\newtheorem{def-lem}[subsection]{Lemma-Definition}

\theoremstyle{remark}

\theoremstyle{plain}

\numberwithin{equation}{subsection}

\def\boxit#1#2{\setbox1=\hbox{\kern#1{#2}\kern#1}%
\dimen1=\ht1 \advance\dimen1 by #1 \dimen2=\dp1 \advance\dimen2 by
#1
\setbox1=\hbox{\vrule height\dimen1 depth\dimen2\box1\vrule}%
\setbox1=\vbox{\hrule\box1\hrule}%
\advance\dimen1 by .4pt \ht1=\dimen1 \advance\dimen2 by .4pt
\dp1=\dimen2 \box1\relax}

\renewcommand{\theequation}{\thesubsection.\arabic{equation}}

\DeclareMathOperator*{\grad}{grad}

\def\ord{{\rm ord}}

\begin{document}

\author[Raf Cluckers]{Raf Cluckers}
\address{Katholieke Universiteit Leuven,
Departement wiskunde, Celestijnenlaan 200B, B-3001 Leu\-ven,
Bel\-gium. Current address: \'Ecole Normale Sup\'erieure,
D\'epartement de ma\-th\'e\-ma\-ti\-ques et applications, 45 rue
d'Ulm, 75230 Paris Cedex 05, France}
\email{raf.cluckers@wis.kuleuven.be}
\urladdr{www.dma.ens.fr/$\sim$cluckers/}

\thanks{$^{\mathrm{*)}}$During the realization of this project, the author was
a postdoctoral fellow of the Fund for Scientific Research - Flanders
(Belgium) (F.W.O.) and was supported by The European Commission -
Marie Curie European Individual Fellowship with contract number HPMF
CT 2005-007121.}

\begin{abstract}
We prove the intersection of Igusa's Conjecture of [Igusa, J.,
\textit{Lectures on forms of higher degree}, Lect.~math.~phys.,
Springer-Verlag, {\bf{59}} (1978)] and the Denef - Sperber
Conjecture of [Denef, J. and Sperber, S., \textit{Exponential sums
mod $p^n$ and {N}ewton polyhedra}, Bull.~Belg.~Math.~Soc.,
{\bf{suppl.}} (2001) 55-63] on nondegenerate exponential sums modulo
$p^m$.
\end{abstract}


\title[{I}gusa and  {D}enef-{S}perber Conjectures]{{I}gusa and {D}enef-{S}perber
Conjectures on nondegenerate $p$-adic exponential sums}

\maketitle

\section{Introduction}

Let $f$ be a polynomial over $\ZZ$ in $n$ variables. Consider the
``global'' exponential sum
$$
S_f(N):= \frac{1}{N^n}\sum_{x\in \{0,\ldots,N-1\}^n}\exp(2\pi
i\frac{ f(x)}{N}),
$$
where $N$ varies over the positive integers. In order to bound
$|S_f(N)|$ in terms of $N$, it is enough to bound
$$|S_f(p^m)|$$
in terms of $m>0$ and prime numbers $p$. When $f$ is nondegenerate
in several senses related to its Newton polyhedron, specific bounds
which depend uniformly on $m$ and $p$ have been conjectured by Igusa
and by  Denef - Sperber.

We prove these bounds, thus solving a conjecture by Denef and
Sperber from a 1990 manuscript \cite{DenSperman} (published in 2001
\cite{DenSper}), and the nondegenerate case of Igusa's conjecture
for exponential sums from the introduction of his book
\cite{Igusa3}.

One of the main points of this article is that, while for finite
field exponential sums like $S_f(p)$ one knows that the weights and
Betti numbers have some uniform behaviour for big $p$, for $p$-adic
exponential sums $S_f(p^m)$ one does not yet completely know what
the analogues of the weights and Betti numbers are, let alone that
they have some uniform behaviour in $p$.

\subsection*{Notation}

Let $f$ be a nonconstant polynomial over $\ZZ$ in $n$ variables with
$f(0)=0$.\footnote{When $f(0)\not=0$, then there is no harm in
replacing $f$ by $f-f(0)$: all corresponding changes in the paper
are easily made.} Write $f(x)=\sum_{i\in \NN^n} a_i x^i$ with
$a_i\in \ZZ$.
 The \emph{global Newton polyhedron} $\Delta_{\rm global}(f)$ of $f$
is the convex hull of the support ${\rm Supp}(f)$ of $f$, with
$${\rm Supp}(f):=\{i\mid i\in\NN^n,\, a_i\not=0\}.$$
  The \emph{Newton polyhedron $\Delta_0(f)$
of $f$ at the origin} is 
$$
\Delta_0(f):=\Delta_{\rm global}(f)+\RR_+^n
$$
with $\RR_+=\{x\in\RR\mid x\geq 0\}$ and $A+B=\{a+b\mid a\in A,\
b\in B\}$ for $A,B\subset\RR^n$.
 For a subset $I$ of $\RR^n$ define
 $$f_I(x):= \sum_{i\in I\cap
\NN^n} a_i x^i.$$
  By the \emph{faces} of $I$ we mean
$I$ itself and each nonempty convex set of the form
$$
\{x\in I\mid  L(x)=0\}
$$
where $L(x)=a_0+ \sum_{i=1}^n a_ix_i$ with $a_i\in \RR$ is such that
 $L(x)\geq 0$ for each $x\in I$.
 By \emph{the proper faces} of $I$ we mean the faces of $I$ that are
different from $I$.
 For $\cI$ a collection of subsets of $\RR^n$, call $f$
\emph{nondegenerate with respect to $\cI$} when $f_I$ has no
critical points on $(\CC^\times)^n$ for each $I$ in $\cI$, where
$\CC^\times=\CC\setminus\{0\}$.
 For $k\in\RR^n_+$ put 
\begin{eqnarray*}
\nu(k) & = & k_1+k_2+\ldots+k_n,\\
N(f)(k) & = & \min_{i\in\Delta_0(f)} k\cdot i,\\
F(f)(k) &= & \{i\in \Delta_0(f)\mid k\cdot i = N(f)(k)\},
\end{eqnarray*}
where $k\cdot i$ is the standard inproduct on $\RR^n$.
Denote by $F_0(f)$ the smallest face of $\Delta_0(f)$ which has
nonempty intersection with the diagonal $\{(t,\ldots,t)\mid
t\in\RR\}$ and let $(1/\sigma(f),\ldots,1/\sigma(f))$ be the
intersection point. If there is no confusion, we write $\sigma$
instead of $\sigma(f)$, $N(k)$ instead of $N(f)(k)$, and $F(k)$ for
$F(f)(k)$. Let $\kappa$ be the codimension of $F_0(f)$ in $\RR^n$.


\section{The main results}
From here up to section \ref{ana}, let $f$ be a nonconstant
polynomial over $\ZZ$ in $n$ variables with $f(0)=0$ (the adaptation
to the case $f(0)\not=0$ is easy).
\begin{theorem}\label{1.2}
Suppose that $f$ is homogeneous and nondegenerate w.r.t. the 
faces of $\Delta_0(f)$. Then there exists $c>0$ and $M>0$ such that
$$
|S_{f}(p^m)|\leq c\, p^{-\sigma m}\, m^{\kappa-1}
$$
for all primes $p>M$ and all integers $m>0$, with $\sigma=\sigma(f)$
and $\kappa$ as in the section on notation. Moreover, $c$ can be
taken depending on $\Delta_0(f)$ only.
\end{theorem}

One sees that the dependence on $p$ and $m$ is very simple. Since
moreover for each $p$ there exists $c_p>0$ such that for all $m>0$
$$
|S_{f}(p^m)|\leq c_pp^{-\sigma m}m^{n-1},
$$
by \cite{Igusa3}, \cite{DenefVeys}, or \cite{DenefBour}, and
properties of toric resolutions, and since $\kappa\leq n$  one finds
the following.

\begin{cor}\label{corig}With $f$ as in Theorem \ref{1.2}
there exists $c>0$ such that for all primes $p$ and all integers
$m>0$,
$$
|S_{f}(p^m)|\leq cp^{-\sigma m}m^{n-1}.
$$
\end{cor}
Denef and Sperber \cite{DenSper} prove Theorem \ref{1.2} under the
extra condition that no vertex of $F_0(f)$ belongs to $\{0,1\}^n$.
Corollary \ref{corig} is the nondegenerate case of Igusa's
conjecture \cite{Igusa3} on exponential sums for toric resolutions
of $f$. The exponent $\sigma$  in the bounds of Theorem \ref{1.2} is
conjectured to be optimal for infinitely many $p$ and $m$ by Denef
and Sargos \cite{DenSar1}, \cite{DenSar2}, in analogy to conjectures
on the real case. When no vertex of $F_0(f)$ belongs to
$\{0,1,2\}^n$, the bounds in Theorem \ref{1.2} are shown to be
optimal   for infinitely many $p$ and $m$ in \cite{DenSper}, Theorem
(1.3).

\section{Denef - Sperber Formula for $S_{f}(p^m)$ for big $p$} 

The following proposition has the same proof as Proposition (2.1) of
\cite{DenSper}, but is slightly more general. We give the proof for
the convenience of the reader.

\begin{prop}\label{r21}Suppose that $f$ is
nondegenerate w.r.t. (all) the faces of $\Delta_0(f)$. Then there
exists $M>0$ such that
 \begin{equation}\label{DF}
S_{f}(p^m)=(1-p^{-1})^n\sum_{\small \mbox{$\tau$ face of
$\Delta_0(f)$}}\big( A(p,m,\tau)
 +
 E(p,f_\tau) B(p,m,\tau)
 \big)
\end{equation}
for all primes $p>M$ and all integers $m>0$, with
$$
A(p,m,\tau):=\sum_{\small
\begin{array}{c}k\in \NN^n
\\ F(f)(k)=\tau\\ N(f)(k)\geq m \end{array}}
p^{-\nu(k)},
$$
$$
B(p,m,\tau):= \sum_{\small
\begin{array}{c}k\in \NN^n
\\ F(f)(k)=\tau\\ N(f)(k)= m-1 \end{array}}p^{-\nu(k)},
$$
and
\begin{equation}\label{Ept}
E(p,f_\tau):=\frac{1}{(p-1)^n} \sum_{\small \mbox{$x\in
\{1,\ldots,p-1\}^n$}}\exp\big(\frac{2\pi i}{p}f_\tau(x)\big).
\end{equation}
\end{prop}
\begin{proof}
Writing
$$
S_f(p^m) = \int_{\ZZ_p}\exp\big(\frac{2\pi i}{p^m}f(x)\big)|dx|,
$$
with $|dx|$ the normalized Haar measure on $\QQ_p^n$, we deduce
$$
S_f(p^m) = \sum_{\small \mbox{$\tau$ face of
$\Delta_0(f)$}}\sum_{\small
\begin{array}{c}k\in \NN^n
\\ F(f)(k)=\tau\end{array}}
 \int_{\ord\, x=k}\exp\big(\frac{2\pi i}{p^m}f(x)\big)|dx|.
$$
Put $x_j=p^{k_j}u_j$ for $k\in\NN^n$. Then $|dx|=p^{-\nu(k)}|du|$
and
$$
f(x)=p^{N(k)}\big(f_{F(k)}(u)+ p(...)\big),
$$
where the dots take values in $\ZZ_p$ 
and where $N(k)=N(f)(k)$. Hence,
\begin{equation}\label{ein}
S_f(p^m) = \sum_{\small \mbox{$\tau$ face of
$\Delta_0(f)$}}\sum_{\small
\begin{array}{c}k\in \NN^n
\\ F(f)(k)=\tau\end{array}}p^{-\nu(k)} \int_{u\in(\ZZ_p^\times)^n}
 \exp\big(\frac{2\pi i}{p^{m-N(k)}}(f_\tau(u)+ p(...))\big)|du|,
\end{equation}
where $\ZZ_p^\times$ denotes the group of $p$-adic units. Because of
the nondegeneratedness assumptions, for $\tau$ a face of
$\Delta_0(f)$ and $p$ a big enough prime, the reduction $f_\tau\bmod
p$ has no critical points on $(\FF_p^\times)^n$.
Hence, the integral in (\ref{ein}) is zero whenever $m-N(f)(k)\geq
2$.
 When $m-N(k)\leq 0$, the integral over $(\ZZ_p^\times)^n$ in
(\ref{ein}) is just the measure of $(\ZZ_p^\times)^n$ and thus
equals $(1-p^{-1})^n$. When $m-N(k)=1$ the integral over
$(\ZZ_p^\times)^n$ in (\ref{ein}) equals $p^{-n}(p-1)^nE(p,f_\tau)$.
 The Proposition now follows from (\ref{ein}).
\end{proof}

\section{Lower bounds for $\nu(k)$}

The main result of this section is:
\begin{theorem}\label{len}
Let $\tau$ be a face of $\Delta_0(f)$. Then one has for all $k$ with
$F(f)(k)=\tau$ that
\begin{equation}\label{enu0}
\nu(k)\geq \sigma(f)\big(N(f)(k) + 1\big) - \sigma(f_\tau),
\end{equation}
where $\sigma(f)$ and $\sigma(f_\tau)$, as well as $\nu(k)$ and
$N(f)(k)$ are as in the section on notation.
\end{theorem}
The main points are that one subtracts $\sigma(f_\tau)$ instead of
$\sigma(f)$ and that $\sigma(f_\tau)\leq \sigma(f)$. Subtracting
$\sigma(f)$ would yield trivial bounds since one has $\nu(k)\geq
\sigma(f) N(k)$ for all $k\in \RR^n_+$. The Theorem's proof is based
on two facts:
\begin{lem}\label{lg}
Let $\tau$ be a face of $\Delta_0(f)$, and let $R_j\in\RR^n$ be
finitely many points belonging to $\tau$. Let $\beta_j\geq 0$
satisfy
$$
\sum_j\beta_j R_j\leq (1/\sigma,\ldots,1/\sigma),
$$
where $a\leq b$ for $a,b\in\RR^n$ means $a_i\leq b_i$ for all $i$.
Then
$$
\sum\beta_j\leq 1.
$$
\end{lem}

\begin{proof}
Clearly there is no point $S$ in the interior of $\Delta_0(f)$ that
satisfies $S\leq (1/\sigma,\ldots,1/\sigma)$. When
$\sum_j\beta_j>1$, then $\sum_j\beta_j R_j$ lies in the interior of
$\Delta_0(f)$.
\end{proof}

\begin{cor}\label{cg}
Let $\tau$, $R_j$, and $\beta_j$ be as in Lemma \ref{lg}. Then
\begin{equation}\label{sst}
\sum\beta_j\leq \frac{\sigma(f_\tau)}{\sigma}.
\end{equation}
\end{cor}
\begin{proof}
Since $ \sum_j\beta_j R_j\leq (1/\sigma,\ldots,1/\sigma)$ one has
$$
\frac{\sigma}{\sigma(f_\tau)}\sum_j\beta_j R_j\leq
(1/\sigma(f_\tau),\ldots,1/\sigma(f_\tau)).
$$
 Lemma \ref{lg} thus implies
$ \frac{\sigma}{\sigma(f_\tau)}\sum_j\beta_j\leq 1. $
\end{proof}

\begin{proof}[Proof of Theorem \ref{len}]
Since $(1/\sigma,\ldots,1/\sigma)$ lies in the interior of $F_0(f)$,
by convexity one can write
$$
(1/\sigma,\ldots,1/\sigma)=\sum_i\alpha_i P_i +\sum_j\beta_jR_j
$$
for some $\alpha_i\geq 0$ and $\beta_j\geq 0$ with
$\sum_i\alpha_i+\sum_j\beta_j =1$ and with $P_i$ finitely many
integral points of $F_0(f)\setminus \tau$ and $R_j$ finitely many
integral points of $\tau$. For $k\in\NN^k$ with $F(f)(k)=\tau$
calculate
\begin{eqnarray}
\nu(k)
 & = & \sigma (1/\sigma,\ldots,1/\sigma)\cdot k\\
 & = &  \sigma \big(\sum_i\alpha_i P_i +\sum_j\beta_jR_j\big)\cdot k\\
 & = & \sigma \big(\sum_i\alpha_i P_i\cdot k +\sum_j\beta_jR_j\cdot k\big)\\
 & \geq & \sigma \big(\sum_i\alpha_i (N(k)+1) +\sum_j\beta_jN(k)\big)\label{eee}\\
 & = & \sigma \big((\sum_i\alpha_i+\sum_j\beta_j) (N(k)+1) -\sum_j\beta_j\big)\\
 & = & \sigma \big((N(k)+1) -\sum_j\beta_j\big)\\
 & \geq & \sigma \big((N(k)+1) -
 \frac{\sigma(f_\tau)}{\sigma}\big)\label{eeee}
\end{eqnarray}
where (\ref{eee}) follows from $k\cdot R_j=N(k)$ and $k\cdot P_i\geq
N(k)+1$ which is true by definition of $N(k)$, and where
(\ref{eeee}) follows from Corollary \ref{cg}.
\end{proof}

\section{Upper bounds for $A(p,m,\tau)$ and $B(p,m,\tau)$}

We recall one result from \cite{DenSper}.
\begin{lem}[\cite{DenSper}, Lemma (3.3)]\label{l3.3}
Let $C$ be a convex polyhedral cone in $\RR_+^n$ generated by
vectors in $\NN^n$, and let $L$ be a linear form in $n$ variables
with coefficients in $\NN$. We denote by $C^{\rm int}$ the interior
of $C$ in the sense of Newton polyhedra. Let $\sigma>0$ and
$\gamma\geq 0$ be real numbers satisfying
\begin{equation}
\nu(k)\geq L(k)\sigma + \gamma,\ \mbox{ for all }k\in C^{\rm
int}\cap \NN^n.
\end{equation}
Put
$$
e=\dim \{k\in C \mid \nu(k)=L(k)\sigma\}.
$$
Then there exists a real number $c>0$ such that for all $m\in\NN$
and for all $p\in\RR$, with $p\geq 2$,
\begin{equation}
\sum_{\begin{array}{c}k\in C^{\rm int}\cap\NN^n
\\ L(k)=m\end{array}} p^{-\nu(k)}\leq
cp^{-m\sigma-\gamma}(m+1)^{\max(0,e-1)}.
\end{equation}
\end{lem}
From this Lemma and from Theorem \ref{len} follows:
\begin{cor}\label{AB}
Let $f$, $A(p,m,\tau)$, and $B(p,m,\tau)$ be as in Proposition
\ref{r21}. Then there exists a real number $c>0$ such that for all
integers $m>0$, for all faces $\tau$ of $\Delta_0(f)$, and for all
big enough primes $p$
\begin{equation}\label{eAe}
A(p,m,\tau)\leq cp^{-m\sigma}m^{\kappa-1}
\end{equation}
 and
\begin{equation}\label{eBe}
B(p,m,\tau)\leq cp^{-m\sigma + \sigma(f_\tau)}m^{\kappa-1}.
\end{equation}
\end{cor}
\begin{proof}
To derive (\ref{eAe}) from Lemma \ref{l3.3}, note that $\nu(k)\geq
N(k)\sigma$ for any $k\in\NN^n$, and that $\kappa= \dim \{k\in
\RR_+^n \mid \nu(k)=N(k)\sigma\}\geq 1$.

To derive (\ref{eBe}) from Lemma \ref{l3.3} and Theorem \ref{len},
use for $C$ the topological closure of the convex hull of $\{0\}\cup
\{k\in\NN^n\mid F(k)= \tau\}$, and note that $C^{\rm int}\cap\NN^n=
\{k\in\NN^n\mid F(k)= \tau\}$. Clearly $\kappa\geq 1$ and
$\kappa\geq \dim \{k\in C \mid \nu(k)=N(k)\sigma\}$. By
(\ref{enu0}), $\nu(k)\geq N(k)\sigma+ \sigma - \sigma(f_\tau)$ for
all $k\in C^{\rm int}\cap\NN^n$.

\end{proof}

\section{Upper bounds for $\sigma(f)$ and $E(p,f_\tau)$}

By Theorem 4 of Katz \cite{Katz}, for $f(x)$ a nonconstant
homogeneous polynomial in $n$ variables over $\ZZ$, and $d$ the
dimension of $\grad f=0$ in $\AA_\CC^n$, there exists $c$ such that
for all big enough $p$ one has
\begin{equation}\label{Ka}
 |\sum_{\small \mbox{$x\in \AA^n(\FF_p)$}}\exp\big(\frac{2\pi
i}{p}f(x)\big)| \leq c p^{\frac{n+d}{2}}.
 \end{equation}
Moreover, $c$ can be taken depending on the degree of $f$ only. This
implies:
\begin{cor}\label{cK}
Suppose that $f(x)$ is homogeneous of degree $\geq 2$ and let $d$ be
the dimension of $\grad f=0$ in $\AA_\CC^n$. Then there exists $c$
such that for all $p$ big enough
\begin{equation}\label{ep}
|\sum_{\small \mbox{$x\in \GG_m^n(\FF_p)$}}\exp\big(\frac{2\pi
i}{p}f(x)\big)| <cp^{\frac{n+d}{2}}
\end{equation}
and hence, for some $c'$ one has, for all  big enough $p$,
\begin{equation}\label{eEp}
|E(p,f)| <c'p^{\frac{-n+d}{2}}
\end{equation}
with $E(p,f)$ as defined by (\ref{Ept}). Moreover, $c$ and $c'$ can
be taken depending on $\Delta_0(f)$ only.
\end{cor}
\begin{proof}
Let $f_0(x_2,\ldots,x_n)$ be the polynomial $f(0,x_2,\ldots,x_n)$.
Clearly $f_0$ is homogeneous in $n-1$ variables. By Katz' result
(\ref{Ka}) it is enough to show that $n-1+d(f_0)\leq n+d$, with
$d(f_0)$ the dimension of $\grad f_0=0$ in $\AA_\CC^{n-1}$. This
inequality follows from writing
$$
f(x)=x_1g(x)+f_0(x_2,\ldots,x_n)
$$
with $g$ a polynomial in $x$, and comparing $\grad f$ with $\grad
f_0$.
\end{proof}

\subsection{}
Let $\{(N_i,\nu_i)\}_{i\in I}$ be the numerical data of a resolution
$h$ of $f$ with normal crossings (that is, if $\pi_{f}:Y\to
\AA^n_\CC$ is an embedded resolution of singularities with normal
crossings of $f=0$, then, for each irreducible component $E_i$ of
$\pi_f^{-1}\circ f^{-1}(0)$, $i\in I$, let $N_i$ be the multiplicity
of $E_i$ in ${\rm div}( f\circ \pi_f)$, and $\nu_i-1$ the
multiplicity of $E_i$ in the divisor associated to
$\pi_f^*(dx_1\wedge...\wedge dx_n)$, cf.~\cite{DenefBour}).  The
\emph{essential numerical data} of $\pi_f$ are the pairs
$(N_i,\nu_i)$ for $i\in J$ with $J=I\setminus I'$ and where $I'$ is
the set of indices $i$ in $I$ such that $(N_{i},\nu_{i})=(1,1)$ and
such that $E_i$ does not intersect another $E_j$ with
$(N_{j},\nu_{j})=(1,1)$.
Define $\alpha(\pi_f)$ as
\begin{equation}\label{alphapi}
\alpha(\pi_f)=-\min_{i\in J}\frac{\nu_i}{N_i}
\end{equation}
when $J$ is nonempty and define $\alpha(\pi_f)$ as $-2n$ otherwise.

It follows from \cite{Cigumodp}, Theorem 5.1 and Corollary 3.4, that
\begin{equation}\label{Ce}
\alpha(\pi_f) \geq \frac{-n+d}{2},
\end{equation}
with $d$ the dimension of $\grad f=0$ in $\AA_\CC^n$, and where the
empty scheme has dimension $-\infty$.


\begin{lem}\label{Epf0} Let $f$ be homogeneous of degree $\geq 2$ and
nondegenerate w.r.t.
the 
faces of $\Delta_0(f)$. Let $d$ be the dimension of $\grad f=0$ in
$\AA_\CC^n$. 
Then
\begin{equation}\label{eps}
\sigma(f)\leq \frac{n-d}{2}.
\end{equation}
\end{lem}
\begin{proof}
By properties of a toric resolution $\pi_f$ of $f=0$, one has that
$$
\sigma(f)=-\alpha(\pi_f),
$$
with $\alpha(\pi_f)$ as defined by (\ref{alphapi}). Now use
(\ref{Ce}).

\end{proof}

From (\ref{eps}) and Corollary \ref{cK}, applied to $f_\tau$,
follows:
\begin{cor}\label{Epf1}
Let $f$ be a homogeneous polynomial of degree $\geq 2$ which is
nondegenerate w.r.t.
the 
faces of $\Delta_0(f)$.
Then there exists $c$ such that for all faces $\tau$ of
$\Delta_0(f)$ and all big enough primes $p$
\begin{equation}\label{Esig}
|E(p,f_\tau)| <cp^{-\sigma(f_\tau)}
\end{equation}
with $E(p,f_\tau)$ as defined by (\ref{Ept}). Moreover, $c$ can be
taken depending on $\Delta_0(f)$ only.
\end{cor}

\section{Proof of the main theorem}

\begin{proof}[Proof of Theorem \ref{1.2} ]
When the degree of $f$ is $\geq 2$, use Proposition \ref{r21},
Corollary \ref{AB}, and (\ref{Esig}). For linear $f$ the theorem is
trivial.
\end{proof}


\section{Comparison with the Denef-Sperber approach}

As mentioned above, Denef and Sperber \cite{DenSper} prove Theorem
\ref{1.2} under the extra condition that no vertex of $F_0(f)$
belongs to $\{0,1\}^n$. Key points in our proof of Theorem \ref{1.2}
are (\ref{enu0}) (which implies Corollary \ref{AB}) and
(\ref{Esig}). Instead of (\ref{enu0}), Denef and Sperber used their
result that, for similar $k$ as in (\ref{enu0}) but assuming the
extra condition that no vertex of $F_0(f)$ belongs to $\{0,1\}^n$,
\begin{equation}\label{DSn}
\nu(k)\geq \sigma(f)(N(f)(k)+1) - \frac{\dim \tau+1}{2}.
\end{equation}
This often fails if one omits the extra condition, see Examples (1)
and (2) below. Instead of (\ref{Esig}), they used the
Adolphson-Sperber \cite{AdoSper}, Denef-Loeser \cite{DenLoexp}
bounds
\begin{equation}\label{DSE}
|E(p,f_\tau)| <cp^{\frac{-\dim \tau-1}{2}}
\end{equation}
which hold (in particular) under the same conditions as for
(\ref{Esig}), but which are sometimes not as good as the bounds
(\ref{Esig}).\footnote{Although (\ref{DSE}) is sometimes sharper
than (\ref{Esig}) in cases where it does not matter for our course.}

 We give two examples where
our methods really make a difference with (\ref{DSn}) and
(\ref{DSE}).\\

\textbf{Examples.}
\begin{itemize}
\item[(1)] First, for $f(x,y,z,u)=xy+zu$ and $\tau=F_0(f)$, one has $\dim
\tau=1$, $\sigma(f_\tau)=\sigma=2$,  (\ref{DSn}) does not hold and
(\ref{DSE}) is not optimal, while (\ref{Esig}) yields the optimal
$|E(p,f_\tau)| <cp^{-2}.$

\item[(2)] Secondly, for $f(x,y,z,u)=xy+zu+xz+ayu$ with $a\in\ZZ$,
$a\not=1$, and $\tau=F_0(f)$, one has $\dim \tau=2$,
$\sigma(f_\tau)=\sigma=2$, (\ref{DSn}) does not hold and (\ref{DSE})
is not optimal, while (\ref{Esig}) yields again the optimal
$|E(p,f_\tau)| <cp^{-2}$ for big $p$. In this example, $E(p,f_\tau)$
can be calculated by performing a transformation on $\GG_m^4$ coming
from an element of $GL_n(\ZZ)$ transforming $f(x)$ into
$f(x',y',z',u')=x'+y'+z'+ax'y'z'^{-1}$; the bounds for $E(p,f_\tau)$
are surprisingly sharp compared, for example, to bounds for the
resembling Kloosterman sums.
\end{itemize}

\section{Analogues over finite extensions of $\QQ_p$ and over
$\FF_q((t))$}\label{ana}

For any nonarchimedean local field $K$ with valuation ring $\cO_K$,
write $\psi_K$ for an additive character
$$
\psi_K:K\to\CC^\times
$$
that is trivial on $\cO_K$ but nontrivial on some element of $K$ of
order $-1$. Write $\ord_K:K^\times\to\ZZ$ for the valuation,
$|\cdot|_K:K\to\RR$ for the norm on $K$, and $\bar K$ for its
residue field, with $q_K$ elements.
 Let $k$ be a number field with ring of integers $\cO_k$. In this
section $f$ is a nonconstant polynomial over $\cO_k[1/N]$ in $n$
variables, with $f(0)=0$ and $N\in\ZZ$.
 For $K$ any nonarchimedean local field that is an algebra over $\cO_k[1/N]$
and for $y\in K^\times$, consider the exponential integral
$$
S_{f,K}(y):= \int_{\cO_K^n}\psi_K(yf(x))|dx|_K,
$$
with $|dx|_K$ the normalized Haar measure on $K^n$. Note that $K$
may be of positive characteristic.
 Then the following generalization of Theorem \ref{1.2}
holds:
\begin{theorem}\label{1.2bbis}
Suppose that $f$ is a homogeneous polynomial over $\cO_k[1/N]$ which
is
nondegenerate w.r.t. the 
faces of $\Delta_0(f)$. Then there exist $c>0$ and $M>N$ such that
$$
|S_{f,K}(y)|\leq c\, |y|_K^{-\sigma}\, |\ord_K(y)|^{\kappa-1}
$$
for all nonarchimedean local fields $K$ that are algebras over
$\cO_k[1/N]$ and have residue characteristic $>M$, and all $y\in
K^\times$ with $\ord_K(y)<0$, with $|\cdot|$ the complex norm.
Moreover, $c$ can be taken depending on $\Delta_0(f)$ only.
\end{theorem}
\begin{proof}
Same proof as of Theorem \ref{1.2}, using Proposition \ref{r21b}
instead of Proposition \ref{r21}.
\end{proof}

\begin{prop}\label{r21b}
Suppose that $f$ is nondegenerate w.r.t. (all) the faces of
$\Delta_0(f)$. Then there exists $M>N$ such that
$$
S_{f,K}(y)=(1-q_K^{-1})^n \sum_{\small \mbox{$\tau$ face of
$\Delta_0(f)$}}\big(A(q_K,m,\tau)
 +
 E(\bar K,\tau,y)B(q_K,m,\tau)
 \big)
$$
for all nonarchimedean local fields $K$ that are algebras over
$\cO_k[1/N]$ and have residue characteristic $>M$ and all $y\in
K^\times$ with $\ord_K(y)\leq 0$.\footnote{For such $K$ and for
$y\in K^\times$ with $\ord_K(y)\geq 0$, one has $ S_{f,K}(y)=1$.}

\par
 In these formulas, $A(q_K,m,\tau)$ and $B(q_K,m,\tau)$ are as in
Proposition \ref{r21}, and
\begin{equation}\label{eEg}
E(\bar K,\tau,y) := \frac{1}{(q_K-1)^n} \sum_{\small \mbox{$u\in
\GG_m^n(\bar K)$}}\psi_{y}(f_\tau(u)),
\end{equation}
with $\psi_{y}$ a nontrivial additive character on $\bar K$
depending on $y$ and $\psi_K$.\footnote{In fact, the character
$\psi_y$ only depends on $\psi_K$ and on $\ac(y)$ for any
multiplicative homomorphism $\ac: K^\times\to \bar K^\times$
extending the natural projection $\cO_K^\times\to \bar K^\times$.}
\end{prop}
\begin{proof}
Same proof as of Proposition \ref{r21}.
\end{proof}

\subsection*{Acknowledgment}
I would like to thank J.~Denef, E.~Hrushovski, and F.~Loeser for
inspiring discussions during the preparation of this paper.

\bibliographystyle{amsplain}
\bibliography{anbib}
\end{document}